\documentclass{article}
\usepackage{amsmath}

\newtheorem{theorem}{Theorem}

\newenvironment{proof}[1][Proof]{\noindent\textbf{#1.} }{\ \rule{0.5em}{0.5em}}
\input{tcilatex}
\begin{document}

\title{{\Large \textbf{The reduction of Laplace equation in certain
Riemannian spaces and the resulting Type II hidden symmetries}}}
\author{Andronikos Paliathanasis\thanks{%
Email: anpaliat@phys.uoa.gr} \ and Michael Tsamparlis\thanks{%
Email: mtsampa@phys.uoa.gr} \\
{\small \textit{Faculty of Physics, Department of Astrophysics - Astronomy -
Mechanics,}}\\
{\small \textit{\ University of Athens, Panepistemiopolis, Athens 157 83,
GREECE}}}
\date{}
\maketitle

\begin{abstract}
We prove a general theorem which allows the determination of Lie symmetries
of Laplace equation in a general Riemannian space using the conformal group
of the space. Algebraic computing is not necessary. We apply the theorem in
the study of the reduction of Laplace equation in certain classes of
Riemannian spaces which admit a gradient Killing vector, a gradient
Homothetic vector and a special Conformal Killing vector. In each reduction
we identify the source of Type II hidden symmetries. We find that in general
the Type II\ hidden symmetries of Laplace equation are directly related to
the transition of the CKVs from the space where the original equation is
defined to the space where the reduced equation resides. In particular we
consider the reduction of Laplace equation (i.e. the wave equation) in
Minkowski space and obtain the results of all previous studies in a
straightforward manner. We consider the reduction of Laplace equation in
spaces which admit Lie point symmetries generated from a non gradient HV and
a proper CKV and we show that the reduction with these vectors does not
produce Type II hidden symmetries. We apply the results to General
Relativity and consider the reduction of Laplace equation in locally
rotational symmetric space times (LRS) and in algebraically special vacuum
solutions of Einstein's equations which admit a homothetic algebra acting
simply transitively. In each case we determine the Type II hidden symmetries.
\end{abstract}


PACS - numbers: 02.20.Sv,~02.30.Jr, 02.40.Ky

\section{Introduction}

Lie point symmetries provide a systematic method to facilitate the solution
of differential equations (DE) because they provide the first order
invariants which can be used to reduce the DE. The reduction is different
for ordinary differential equations (ODEs) and for partial differential
equations (PDEs). In the case of ODEs the use of a Lie point \ symmetry
reduces the order of the ODE by one while in the case of PDEs reduces by one
the number of independent and dependent variables, but not the order of the
PDE. A common characteristic in the reduction of both cases is that the Lie
point symmetry which is used for the reduction is not admitted as such by
the reduced DE, it is "lost". In this procedure it is possible that the
reduced equation has more Lie point symmetries than the original equation.
These new point symmetries have been called Type II\ hidden symmetries and
have been the subject of numerous papers during the recent years (e.g. \cite%
{AbrahamGuo1994,AbrahamGovinderLeach1995,LeachGovinderAbraham1999,AbrahamGovinder2008,ASG,AbrahamGovinderArrige2006,Abraham2007,Moitsheko(2004)}%
).

The importance of Type II\ hidden symmetries is that they can be used to
reduce further the reduced equation. Concerning the origin of Type II\
hidden symmetries, both for ODEs and PDEs, it has been shown that they can
be viewed as having two sources. Either the point and the nonlocal /
generalized symmetries of \emph{a given} higher order equation, or the point
symmetries of \emph{a variety} of higher order DEs which reduce to this
particular DE (see \cite%
{ASG,LeachGovinderAndriopoulos2012,Govinder2001,GovinderAbraham2009} and
references therein).

In the present paper we study the reduction and the consequent existence of
Type II\ hidden symmetries\ of Laplace equation in certain classes of
Riemannian spaces. In particular we prove a general theorem which allows us
to study the reduction of Laplace equation in a general Riemannian space
using the conformal group of the space. There is no need to employ algebraic
computing and it is enough to work with pure geometric arguments. We note
that the use of algebraic computing for higher dimensions and more complex
metrics can be rather inapplicable, whereas the theorem applies
irrespectively of the (finite) dimension and the complexity of the metric.

The Laplace equation in a general Riemannian space%
\begin{equation}
\frac{1}{\sqrt{\left\vert g\right\vert }}\frac{\partial }{\partial x^{i}}%
\left( \sqrt{\left\vert g\right\vert }g^{ij}\frac{\partial }{\partial x^{j}}%
\right) u\left( x^{k}\right) =0  \label{PE.9}
\end{equation}%
has two Lie point symmetries, the $u\partial _{u}$ and the $b(x^{k})\partial
_{u}$ where $b(x^{k})$ is a solution of Laplace equation. These two
symmetries are not useful for reduction. In order to find `sound' reductions
of Laplace equation we have to consider Riemannian spaces which admit some
type of symmetry(ies) (these symmetries are not Lie point symmetries and are
called collineations). Indeed it has been shown \cite{Bozhkov} that the Lie
point symmetries of Laplace equation in a Riemannian space are generated
from the conformal algebra of the space. Therefore in order to have more Lie
point symmetries of DEs, hence the possibility of the existence of Type II\
hidden symmetries, we have to work in spaces which admit a conformal algebra.

The structure of the paper is as follows. In section \ref{basicproperties}
we give the basic definitions of Lie point symmetries and certain facts
concerning the conformal algebra of a Riemannian space. In section \ref%
{LSCOl} we prove a general theorem concerning the Lie point symmetries of
the Poisson equation in a Riemannian space. We also give some results which
relate the Lie point symmetries of the Poisson, the Laplace and the Klein
Gordon equations with the conformal algebra of the space. In section \ref%
{Reduction of the Laplace in certain Riemannian spaces} we consider the
reduction of Laplace equation in various general classes of spaces which
admit a nonzero conformal algebra. In particular we consider (a)
decomposable spaces - that is Riemannian spaces which admit a gradient
Killing (equivalently constant) vector field (KV) - (b) spaces which admit a
gradient homothetic vector field (HV) and (c) spaces which admit a special
conformal Killing vector field (sp.CKV). In section \ref{Applications}, we
apply the results of section \ref{Reduction of the Laplace in certain
Riemannian spaces} and determine the Type II symmetries of Laplace equation
in four dimensional Minkowski spacetime and we recover and complete well
known results \cite{AbrahamGovinderArrige2006}. In section \ref{ExLRS} we
consider the reduction of Laplace equation in an LRS spacetime, which is an
important class of spacetimes in General Relativity. In order to study the
reduction of Laplace equation by a non-gradient HV and a proper CKV we
consider two further examples. In section \ref{The Heat equation in spaces
which admit a non-gradient HV} we consider the algebraically special vacuum
solutions of Einstein's equations on which a homothetic vector field acts
transitively \cite{Steele1991(b)} and make the reduction using the Lie point
symmetry generated by the nongradient HV. In section \ref{FRWCKV} \ we do
the same in a conformally flat FRW type space which admits a homothetic
vector field. We reemphasize that all results are derived in a purely
geometric manner without the use of a computer package. However they have
been verified with the libraries PDEtools and SADE \cite{PDEtools,SADE} of
Maple\footnote{%
www.maplesoft.com}. Finally in section \ref{Conclusion} we draw our
conclusions.

\section{Lie point symmetries of a PDE and CKVs of a Riemannian space}

\label{basicproperties}

In this section, we give the definition of Lie point symmetries of a DE and
the definition of conformal Killing vector fields (CKVs) of a Riemannian
space.

\subsection{Lie point symmetries of differential equations}

A partial differential equation (PDE) is a function $H=H\left(
x^{i},u,u_{,i},u_{,ij},..\right)$ in the jet space $\bar{B}_{\bar{M}}\left(
x^{i},u,u_{,i},u_{,ij},..\right)$ where $x^{i}$ are the independent
variables and $u^{A}$ are the dependent variables. An infinitesimal point
transformation%
\begin{eqnarray}
\bar{x}^{i} &=&x^{i}+\varepsilon \xi ^{i}\left( x^{k},u\right)  \label{pde1}
\\
\bar{u} &=&\bar{u}+\varepsilon \eta \left( x^{k},u\right)  \label{pde2}
\end{eqnarray}%
is a Lie point symmetry of the PDE $H\left( x^{i},u,u_{,i},u_{,ij},..\right) 
$ with generator%
\begin{equation}
X=\xi ^{i}\left( x^{k},u^{B}\right) \partial _{t}+\eta ^{A}\left(
x^{k},u^{B}\right) \partial _{B}  \label{pde3}
\end{equation}%
where $\varepsilon $ is an infinitesimal parameter, if there exists a
function $\lambda $ such as the following condition holds \cite{Ibragimov,
Stephani}%
\begin{equation}
X^{\left[ n\right] }H\left( x^{i},u,u_{,i},u_{,ij},..\right) =\lambda
H\left( x^{i},u,u_{,i},u_{,ij},..\right) ~, mod H=0  \label{pde4}
\end{equation}%
$X^{\left[ n\right] }$ is the nth prolongation of (\ref{pde3}) defined as
follows,%
\begin{equation}
X^{\left[ n\right] }=X+\eta _{\left[ i\right] }\partial _{u_{i}}+...+\eta _{%
\left[ ij...i_{n}\right] }\partial _{u_{ij...i_{n}}}
\end{equation}%
where%
\begin{equation*}
\eta _{\left[ i\right] }=D_{i}\eta -u_{,j}D_{i}\xi ^{j}
\end{equation*}%
\begin{equation*}
\eta _{\left[ ij..i_{n}\right] }=D_{i_{n}}\eta _{\left[ ij..i_{n-1}\right]
}-u_{ij..k}D_{i_{n}}\xi ^{k}.
\end{equation*}

From condition (\ref{pde4}) one defines the Lagrange system%
\begin{equation*}
\frac{dx^{i}}{\xi ^{i}}=\frac{du}{\eta }=\frac{du_{i}}{\eta _{\left[ i\right]
}}=...=\frac{du_{ij..i_{n}}}{\eta _{\left[ ij...i_{n}\right] }}
\end{equation*}%
whose solution provides the characteristic functions 
\begin{equation*}
W^{\left[ 0\right] }\left( x^{k},u\right) ,~W^{\left[ 1\right] i}\left(
x^{k},u,u_{i}\right) ,...,W^{\left[ n\right] }\left(
x^{k},u,u_{,i},...,u_{ij...i_{n}}\right) .
\end{equation*}%
The solution $W^{\left[ n\right] }$ is called the nth order invariant of the
Lie symmetry vector (\ref{pde3}). These invariants are used in order to
reduce the order of the PDE (for details see e.g. \cite{Stephani}).

\subsection{Collineations of Riemannian spaces}

\label{Col}

In the following $L_{\xi }$ denotes Lie derivative with respect to the
vector field $\xi ^{i}$.

A vector field $\xi ^{i}$ is a CKV of a metric $g_{ij}$ if $L_{\xi
}g_{ij}=2\psi g_{ij}$. If $\psi =0$ then $\xi ^{i}$ is a KV, if $\psi
_{,i}=0 $ then $\xi ^{i}$ is a HV~and if $\psi _{;ij}=0$ then~$\xi ^{i}$ is
a special CKV (sp.CKV) and also $\psi _{;i}$ is a gradient KV or,
equivalently, a constant vector field. A CKV which is neither of the above
cases (i.e. $\psi _{;ij}\neq 0$ ) is called a proper CKV.

Two metrics $g_{ij},~\bar{g}_{ij}$ are conformally related if $\bar{g}%
_{ij}=N^{2}g_{ij}$ where the function ~$N^{2}$ is the conformal factor. If $%
\xi ^{i}$ is a CKV\ of the metric $\bar{g}_{ij}$ so that $L_{\xi }\bar{g}%
_{ij}=2\bar{\psi}\bar{g}_{ij}$ then $\xi ^{i}$ it is also a CKV of the
metric $g_{ij}$, that is, $L_{\xi }g_{ij}=2\psi g_{ij}$ where the conformal
factor 
\begin{equation*}
\psi =\bar{\psi}N^{2}-NN_{,i}\xi ^{i}.
\end{equation*}%
This means that two conformally related metrics have the same CKVs but
different Killing/homothetic/Sp.CKVs. For example a KV of one is not in
general a KV of the other and so on. This is an important observation which
shall be useful in the following sections. In Appendix \ref{CKVsFlat} we
give the vector fields of the conformal algebra of the flat space in
Cartesian coordinates. Details on the CKVs and their geometric properties
can be found e.g. in \cite{Hallbook}

\section{Lie point symmetries of Laplace equation}

\label{LSCOl}

In a general Riemannian space with metric $g_{ij}$ Poisson equation is%
\begin{equation}
\Delta u-f\left( x^{i},u\right) =0  \label{PE.1}
\end{equation}%
where $\Delta $ is the Laplace operator $\Delta =\frac{1}{\sqrt{g}}\frac{%
\partial }{\partial x^{i}}\left( \sqrt{g}g^{ij}\frac{\partial }{\partial
x^{j}}\right) $ and $q=q(t,x^{k},u)$. Equation (\ref{PE.1}) can also be
written 
\begin{equation}
g^{ij}u_{ij}-\Gamma ^{i}u_{i}\ =f\left( x^{i},u\right)  \label{PE.2}
\end{equation}%
where $\Gamma ^{i}=\Gamma _{jk}^{i}g^{jk}$ and $\Gamma _{jk}^{i}$ are the
Christofell Symbols of the metric $g_{ij}$. The Lie symmetries of the
Poisson equation for $f=f\left( u\right) $ have been given in \cite%
{Bozhkov,Ibragimov}. Here we generalize the result of \cite{Bozhkov} by
considering $f=f\left( x^{i},u\right) $ and give a new concise proof.

\begin{theorem}
\label{Theor}The Lie point symmetries of the Poisson equation (\ref{PE.1})$\ 
$are generated from the CKVs of the metric~$g_{ij}$ of the $\Delta -$%
operator as follows

a) For $n>2$,~the Lie point symmetry vector is 
\begin{equation}
X=\xi ^{i}\left( x^{k}\right) \partial _{i}+\left( \frac{2-n}{2}\psi \left(
x^{k}\right) u+a_{0}u+b\left( x^{k}\right) \right) \partial _{u}
\label{PE.7}
\end{equation}%
where $\xi ^{i}\left( x^{k}\right) $ is a CKV with conformal factor $\psi
\left( x^{k}\right) $ and the following condition holds%
\begin{equation}
\frac{2-n}{2}\Delta \psi u+g^{ij}b_{i;j}-\xi ^{k}f_{,k}-\frac{2-n}{2}\psi
uf_{,u}-\frac{n+2}{2}\psi f-bf_{,u}=0  \label{PE.8}
\end{equation}

b) For $n=2$, the Lie point symmetry vector is%
\begin{equation}
X=\xi ^{i}\left( x^{k}\right) \partial _{i}+\left( a_{0}u+b\left(
x^{k}\right) \right) \partial _{u}
\end{equation}%
where $\xi ^{i}$ is a CKV (i.e. $L_{\xi }g_{ij}=2\psi \left( x^{k}\right)
g_{ij}$) and the following condition holds%
\begin{equation}
g^{ij}b_{;ij}-\xi ^{k}f_{,k}-a_{0}uf_{,u}+\left( a_{0}-2\psi \right)
f-bf_{,u}=0
\end{equation}%
where the function $b$ is solution of Laplace equation.

\begin{proof}
In \cite{JGP} it has been shown\footnote{%
The derivation of these conditions is a straightforward calculation using
the Lie symmetry condition. The details of the calculation can be found in 
\cite{JGP}.} that the Lie point symmetry conditions for the PDE of the form%
\begin{equation}
A^{ij}u_{ij}-B^{i}(x,u)u_{i}-f(x,u)=0  \label{PE.3}
\end{equation}%
are as follows: 
\begin{equation}
A^{ij}(a_{ij}u+b_{ij})-(a_{i}u+b_{i})B^{i}-\xi
^{k}f_{,k}-auf_{,u}-bf_{,u}+\lambda f=0  \label{PE.3a}
\end{equation}%
\begin{equation}
A^{ij}\xi _{,ij}^{k}-2A^{ik}a_{,i}+aB^{k}+auB_{,u}^{k}-\xi
_{,i}^{k}B^{i}+\xi ^{i}B_{,i}^{k}-\lambda B^{k}+bB_{,u}^{k}=0  \label{PE.3b}
\end{equation}%
\begin{equation}
L_{\xi ^{i}\partial _{i}}A^{ij}=(\lambda -a)A^{ij}-\eta A^{ij}{}_{,u}
\label{PE.3c}
\end{equation}%
where the generator of the Lie point symmetry is 
\begin{equation}
\mathbf{X}=\xi ^{i}\left( x^{i}\right) \partial _{i}+\left( a\left(
x^{k}\right) u+b\left( x^{i}\right) \right) \partial _{u}.  \label{PE.6}
\end{equation}

For the Poisson equation (\ref{PE.1}) we have $A^{ij}=g^{ij}~$and $%
B^{i}=\Gamma ^{i}.$ Replacing in conditions (\ref{PE.3a})-(\ref{PE.3c}) we
find 
\begin{equation}
g^{ij}(a_{ij}u+b_{ij})-(a_{,i}u+b_{,i})\Gamma ^{i}-\xi
^{k}f_{,k}-auf_{,u}-bf_{,u}+\lambda f=0  \label{WDW.01}
\end{equation}%
\begin{equation}
g^{ij}\xi _{,ij}^{k}-2g^{ik}a_{,i}+a\Gamma ^{k}-\xi _{,i}^{k}\Gamma ^{i}+\xi
^{i}\Gamma _{,i}^{k}-\lambda \Gamma ^{k}=0  \label{WDW.02}
\end{equation}%
\begin{equation}
L_{\xi }g_{ij}=(a-\lambda )g_{ij}.  \label{WDW.03}
\end{equation}%
Equation (\ref{WDW.02}) becomes (see \cite{Bozhkov}) 
\begin{equation}
g^{jk}L_{\xi }\Gamma _{.jk}^{i}=2g^{ik}a_{,i}.  \label{W3b}
\end{equation}%
Equation (\ref{WDW.03}) gives that $\xi ^{i}$ is a CKV.

a) For $n>2$, since $\xi ^{i}$ is a CKV, equation (\ref{W3b}) becomes 
\begin{equation*}
\frac{2-n}{2}(a-\lambda )^{,i}=2a^{,i}\rightarrow \left( a-\lambda \right)
^{i}=\frac{4}{2-n}a^{,i}.
\end{equation*}%
Therefore,~%
\begin{equation*}
\psi =\frac{2}{2-n}a+a_{0}
\end{equation*}%
where $2\psi =\left( a-\lambda \right) $ is the conformal factor of $\xi
^{i} $. Furthermore we have%
\begin{equation*}
\lambda ^{i}=-\frac{\left( n+2\right) }{\left( 2-n\right) }a^{i}.
\end{equation*}%
Finally from (\ref{WDW.01}) we have the constraint condition%
\begin{equation*}
g^{ij}a_{i;j}u+g^{ij}b_{;ij}-\xi ^{k}f_{,k}-auf_{,u}+\lambda f-bf_{,u}=0.
\end{equation*}

b) For $n=2$ we have that $g^{jk}L_{\xi }\Gamma _{.jk}^{i}=0,$ which implies 
$a_{,i}=0\rightarrow a=a_{0}$. Finally we have%
\begin{equation*}
\lambda =a_{0}-2\psi
\end{equation*}%
and (\ref{WDW.01}) becomes%
\begin{equation*}
g^{ij}b_{;ij}-\xi ^{k}f_{,k}-a_{0}uf_{,u}+\left( a_{0}-2\psi \right)
f-bf_{,u}=0
\end{equation*}
\end{proof}
\end{theorem}

There are two important forms of the Poisson equation which are of special
interest. The Laplace equation (\ref{PE.9}) defined by $f\left(
x^{i},u\right) =0$ and the Klein-Gordon equation defined by $f\left(
x^{i},u\right) =V\left( x^{i}\right) u$ 
\begin{equation}
\Delta u-V\left( x^{k}\right) u=0  \label{PE.10}
\end{equation}%
where $V(x^{i})$ is the `potential'.

The Lie symmetries of (\ref{PE.9}) have been given in \cite{Bozhkov,Stephani}%
. For the convenience of the reader and because we shall make use of this
result in the following, we state this result below.

\begin{theorem}
\label{BOZKOV}The Lie point symmetries of Laplace equation are generated
from the CKVs of the metric $g_{ij}$ as follows:

a) for $n>2$, the Lie point symmetry vector is\newline
\begin{equation}
X=\xi ^{i}\left( x^{k}\right) \partial _{i}+\left[ \left( \frac{2-n}{2}\psi
\left( x^{k}\right) +a_{0}\right) u+b\left( x^{k}\right) \right] \partial
_{u}  \label{PE.11}
\end{equation}%
where $\xi ^{i}\left( x^{k}\right) $ is a CKV with conformal factor $\psi
\left( x^{k}\right) $ and the conformal factor $\psi \left( x^{k}\right) $
and the function $b\left( x^{k}\right) $ satisfy Laplace equation.

b) for $n=2$, the Lie point symmetry vector is 
\begin{equation}
X=\xi ^{i}\left( x^{k}\right) \partial _{i}+\left( a_{0}u+b\left(
x^{k}\right) \right) \partial _{u}
\end{equation}%
where $\xi ^{i}\left( x^{k}\right) $ is a CKV with conformal factor $\psi
\left( x^{k}\right) $ and $b\left( x^{k}\right) $ is a solution of Laplace
equation.
\end{theorem}

The Lie point symmetries of the Klein Gordon equation follow from Theorem %
\ref{Theor}.

\begin{theorem}
\label{KG}The Lie point symmetries of the Klein Gordon equation (\ref{PE.10}%
) are generated from the CKVs of the metric~$g_{ij}$ as follows

a) for $n>2$,~the Lie symmetry vector is\newline
\begin{equation}
X=\xi ^{i}\left( x^{k}\right) \partial _{i}+\left( \frac{2-n}{2}\psi \left(
x^{k}\right) u+a_{0}u+b\left( x^{k}\right) \right) \partial _{u}
\label{PE.13}
\end{equation}%
where $\xi ^{i}$ is a CKV with conformal factor $\psi \left( x^{k}\right) ,$ 
$b\left( x^{k}\right) $ is a solution of the Klein Gordon equation (\ref%
{PE.10}) and the following condition is satisfied%
\begin{equation}
\xi ^{k}V_{,k}+2\psi V-\frac{2-n}{2}\Delta \psi =0  \label{PE.14}
\end{equation}%
b) for $n=2$, the Lie symmetry vector is 
\begin{equation}
X=\xi ^{i}\left( x^{k}\right) \partial _{i}+\left( a_{0}u+b\left(
x^{k}\right) \right) \partial _{u}
\end{equation}%
where $\xi ^{i}$ is a CKV with conformal factor $\psi \left( x^{k}\right) ,$ 
$b\left( x^{k}\right) $ is a solution of the Klein Gordon equation (\ref%
{PE.10}) and the following condition is satisfied%
\begin{equation}
\xi ^{k}V_{,k}+2\psi V=0.
\end{equation}
\end{theorem}

We recall that the conformal algebra of a two-dimensional space has infinite
dimension, a result which we shall need below. Because all two dimensional
spaces are conformally flat they all have the same conformal algebra. In the
subsequent sections we apply theorem \ref{Theor} to study the reduction of
Laplace equation in certain classes of Riemannian spaces which admit a
conformal algebra.

\subsection{The conformal algebra of certain classes of Riemannian spaces}

\label{Laplace equation in certain Riemanian spaces}

In a general Riemannian space Laplace equation (\ref{PE.9}) admits the Lie
point symmetries 
\begin{equation*}
~X_{u}=u\partial _{u}~,~X_{b}=b\left( t,x\right) \partial _{u}
\end{equation*}%
where $b\left( t,x\right) $ is a solution of Laplace equation. These
symmetries are too general to provide useful reductions and lead to reduced
PDEs which posses Type II\ hidden symmetries. However the above theorems
indicate that if we restrict our considerations to spaces which admit a
conformal algebra (proper or not) then we will have new Lie point symmetries
hence new reductions of Laplace equation, which might lead to Type II\
hidden symmetries.

In the following we consider the following classes of Riemannian spaces
which admit a conformal algebra.

a. Spaces which admit a gradient KV (decomposable spaces)

b. Spaces which admit a gradient HV and

c. Spaces which admit a sp.CKV (this case is a subcase of a.).

The generic form of the metric for each type of space is as follows ($%
A,B,\ldots =1,2,\ldots ,n$):

\begin{enumerate}
\item[a.] If an $1+n~$dimensional Riemannian space admits a non null
gradient KV, the $S^{i}=\partial _{z}~\left( S=z\right) $ say, then the
space is decomposable along $\partial _{z}$ and the metric takes the form
(see e.g. \cite{TNA}) 
\begin{equation}
ds^{2}=dz^{2}+h_{AB}y^{A}y^{B}~,~h_{AB}=h_{AB}\left( y^{C}\right)
\label{WH.03}
\end{equation}

\item[b.] If an $n~$dimensional Riemannian space admits a gradient HV, the $%
H^{i}=r\partial _{r}~\left( H=\frac{1}{2}r^{2}\right) ,~\psi _{H}=1$ say,
then the metric can be written in the form \cite{Tupper}%
\begin{equation}
ds^{2}=dr^{2}+r^{2}h_{AB}dy^{A}dy^{B}~,~~h_{AB}=h_{AB}\left( y^{C}\right)
\label{LEH.01}
\end{equation}

\item[c.] If an $n~$dimensional Riemannian space admits a \ non null KV and
a gradient HV, then the space admits a \ sp.CKV and the metric can be
written in the form%
\begin{equation}
ds^{2}=-dz^{2}+dR^{2}+R^{2}f_{AB}\left( y^{C}\right) dy^{A}dy^{B}
\label{LES.01}
\end{equation}%
where the sp.CKV is $C_{S}=\frac{z^{2}+R^{2}}{2}\partial _{z}+zR\partial
_{R} $ with conformal factor $\psi _{C_{S}}=z$ ~\cite{Tupper,Hall}.
\end{enumerate}

The Riemannian spaces which admit a non-gradient proper HV do not have a
generic form for their metric. However the \ spaces for which the HV\ acts
simply transitively are a few and are given together with their homothetic
algebra in \cite{Steele1991(b)}. A\ special class of these spaces are the
algebraically special vacuum solutions of Einstein equations \cite%
{Steele1991(b)}. In the following we shall consider the reduction of Laplace
equation in the\footnote{%
This is not the most general Petrov type III space-time} Petrov type III
solution only with line element%
\begin{equation}
ds^{2}=2d\rho dv+\frac{3}{2}xd\rho ^{2}+\frac{v^{2}}{x^{3}}\left(
dx^{2}+dy^{2}\right)  \label{P3.00}
\end{equation}%
in which the symmetry generated by the non gradient HV~$H=v\partial
_{v}+\rho \partial _{\rho }~,~\psi _{H}=1.$ The reduction of Laplace
equation in the rest of the Petrov types in this class of solutions is
similar both in the working method and the results and there is no need to
consider them explicitly.

Finally we shall consider the conformally flat space 
\begin{equation}
ds^{2}=e^{2t}\left[ dt^{2}-\delta _{AB}y^{A}y^{B}\right]  \label{EC.01}
\end{equation}%
which admits a non-gradient HV.

In what follows all spaces are assumed to be of dimension $n>2$.

\section{Reduction of Laplace equation in certain Riemannian spaces}

\label{Reduction of the Laplace in certain Riemannian spaces}

As we have seen in Theorem \ref{BOZKOV} the Lie point symmetries of Laplace
equation (\ref{PE.9}) in a Riemannian space are the CKVs (not necessarily
proper) whose conformal factor satisfies Laplace equation. This condition is
satisfied trivially by the KVs ($\psi =0),$ the HV ($\psi _{;i}=0)$ and the
sp.CKVs ($\psi _{;ij}=0).$ Therefore these vector fields (which span a
subalgebra of the conformal group) are among the Lie symmetries of Laplace
equation. Concerning the proper CKVs it is not necessary that their
conformal factor satisfies Laplace equation (\ref{PE.9}) therefore they may
not produce Lie point symmetries of Laplace equation.

\subsection{Riemannian spaces admitting a gradient KV}

\label{grKV}

Without loss of generality we assume the gradient KV to be the $\partial
_{z} $ so that the metric has the generic form (\ref{WH.03}) where $h_{AB},$ 
$A,B,C=1,...,n$ is the metric of the $n-$ dimensional space. For the metric (%
\ref{WH.03}) Laplace equation (\ref{PE.9}) takes the form%
\begin{equation}
u_{zz}+h^{AB}u_{AB}-\Gamma ^{A}u_{B}=0  \label{WH.04}
\end{equation}%
and admits as \emph{extra }Lie point symmetry the gradient KV $\partial
_{z}. $

We reduce (\ref{WH.04}) using the zeroth order invariants $y^{A}~,~w=u~$ of
the extra Lie point symmetry $\partial _{z}$. Taking these invariants as new
coordinates eqn (\ref{WH.04}) reduces to 
\begin{equation}
_{h}\Delta w=0  \label{WH.05}
\end{equation}%
which is Laplace equation in the~$n-$dimensional space with metric $h_{AB}.$
Now we recall the (easy to show) result that the conformal algebra of the $n$
metric $h_{AB}$ and the $1+n$ metric (\ref{WH.03}) are related as follows 
\cite{TNA}:

a. The KVs of the $n-$ metric are identical with those of the $n+1$ metric
(apart from the vector $\partial _{z}$).

b. The $1+n$ metric admits a HV\ if and only if the $n$ metric admits one
and if $_{n}H^{A}$ is the HV of the $n$ metric then the HV\ of the $1+n$
metric is given by the expression 
\begin{equation}
_{1+n}H^{\mu }=z\delta _{z}^{\mu }+_{n}H^{A}\delta _{A}^{\mu }\qquad \mu
=x,1,\ldots,n.  \label{WH.05.b}
\end{equation}

c. The \thinspace $1+n$ metric admits CKVs if and only if the $n$ metric $%
h_{AB}$ admits a gradient CKV (for details see proposition 2.1 of \cite{TNA}%
).

Therefore Type II hidden symmetries for (\ref{PE.9}) exist if the $n- $%
metric $h_{AB}$ admits more symmetries. In particular the sp.CKVs of the $%
h_{AB}$ metric as well as the proper CKVs whose conformal function is a
solution of Laplace equation (\ref{WH.05}) generate Type II\ hidden
symmetries.

\subsection{Riemannian spaces admitting a gradient HV}

\label{grHV}

In Riemannian spaces which admit a gradient HV, $H$ say, there exists a
coordinate system in which the metric is written in the form (\ref{LEH.01})
and the gradient HV is $H=r\partial _{r}.$ In these coordinates the
Laplacian takes the form 
\begin{equation}
u_{rr}+\frac{1}{r^{2}}h^{AB}u_{AB}+\frac{\left( n-1\right) }{r}u_{r}-\frac{1%
}{r^{2}}\Gamma ^{A}u_{A}=0  \label{LEH.02}
\end{equation}%
and admits the extra Lie point symmetry $H$ (see Theorem \ref{BOZKOV}). We
reduce (\ref{LEH.02}) using $H.$

The zero order invariants of $H$ are $y^{A}~,~w\left( y^{A}\right) ~$and
using them it follows easily that the reduced equation is%
\begin{equation}
_{h}\Delta w=0  \label{LEH.03}
\end{equation}%
that is, the Laplacian defined by the metric $h_{AB}$.

It is easy to establish the following results concerning the conformal
algebras of the metrics (\ref{LEH.01}) and $h_{AB}$.\newline
1. The KVs of $h_{AB}$ are also KVs of (\ref{LEH.01}).\newline
2. The HV of (\ref{LEH.01}), if it exist, is independent from that of $%
h_{AB} $.\newline
3. The metric (\ref{LEH.01}) admits proper CKVs if and only if the $n$
metric $h_{AB}$ admits gradient CKVs \cite{TNA}. This is because (\ref%
{LEH.01}) is conformally related with the decomposable metric 
\begin{equation}
ds^{2\prime }=d\bar{r}^{2}+h_{AB}\left( y^{C}\right) dy^{A}dy^{B}.
\end{equation}

The above imply that Type II hidden symmetries we shall have from the HV of
the metric $h_{AB},$ the sp.CKVs and finally from the proper CKVs of $h_{AB}$
whose conformal factor is a solution of Laplace equation (\ref{LEH.03}).

\subsection{Riemannian spaces admitting a sp.CKV}

\label{spCKV}

It is known \cite{Hall} that if a decomposable $n=m+1$ dimensional ($n>2)$
Riemannian space which admits the non null gradient KV $K_{G}=\partial
_{z},~ $ also admits sp.CKVs (there is a 1:1 correspondence between the non
null gradient KVs and the sp.CKVs in a Riemannian space see \cite{Hall})
then also admits a gradient HV. In these spaces there exists always a
coordinate system in which the metric is written in the form (\ref{LES.01})
where $\partial _{z}$ is the gradient KV and $z\partial _{z}+R\partial _{R}$
is the gradient HV. $f_{AB}\left( y^{C}\right) $ $A,B,C,..=1,2,,...m-1$ is
an $m-1$ \ dimensional metric. For a general $m-1~$ dimensional metric $%
f_{AB}$ the $n-$dimensional metric (\ref{LES.01}) admits the following
special conformal group 
\begin{eqnarray*}
K_{G} &=&\partial _{z}~,~H=z\partial _{z}+R\partial _{R} \\
C_{S} &=&\frac{z^{2}+R^{2}}{2}\partial _{z}+zR\partial _{R}
\end{eqnarray*}%
where $K_{G}$ is a gradient KV, $H$ is a gradient HV and $C_{S}$ is a sp.CKV
with conformal factor $\psi _{C_{S}}=z$. In these coordinates Laplace
equation takes the form%
\begin{equation}
-u_{zz}+u_{RR}+\frac{1}{R^{2}}h^{AB}u_{AB}+\frac{\left( m-1\right) }{R}u_{R}-%
\frac{1}{R^{2}}\Gamma ^{A}u_{A}=0.  \label{LES.02}
\end{equation}

From Theorem \ref{BOZKOV} we have that the extra Lie point symmetries of (%
\ref{LES.02}) are the vector fields%
\begin{eqnarray}
X^{1} &=&K_{G}~,~X^{2}=H  \label{LES.02c} \\
X^{3} &=&C_{S}+2pzu\partial _{u}
\end{eqnarray}%
where~$2p=\frac{1-m}{2}~.$ The non - zero commutators are%
\begin{equation*}
\left[ X^{1},X^{2}\right] =X^{1}~,~\left[ X^{2},X^{3}\right] =X^{3}~,~\left[
X^{1},X^{3}\right] =X^{2}+2pu\partial _{u}.
\end{equation*}

We consider the reduction of (\ref{LES.02}) with each of the extra Lie
symmetries.

Reduction with the gradient KV $X^{1}$ reduces the Laplacian (\ref{LES.02})
to (\ref{LEH.02}) which admits the Lie point symmetry $X^{2}$ generated by
the HV. This result is expected \cite{AbrahamGovinderArrige2006} because $%
\left[ X^{1},X^{2}\right] =X^{1}$ hence the Lie point symmetry $X^{2}$ is
inherited. Therefore in this reduction the Type II symmetries are generated
from the CKVs of the metric (\ref{LEH.01}). It is possible to continue the
reduction by the gradient HV\ $H$ and then we find the results of section %
\ref{grHV}.

The reduction with a gradient HV has been studied in section \ref{grHV}. To
apply the results of section \ref{grHV} in the present case we have to bring
the metric (\ref{LES.01}) to the form (\ref{LEH.01}). For this we consider
the transformation%
\begin{equation*}
z=r\sinh \theta ~,~R=r\cosh \theta
\end{equation*}%
which brings (\ref{LES.01}) to%
\begin{equation}
ds^{2}=dr^{2}+r^{2}\left( -d\theta ^{2}+\cosh ^{2}\theta
f_{AB}y^{A}y^{B}\right)
\end{equation}%
so that the metric $h_{AB}$ of \ (\ref{LEH.01}) is:%
\begin{equation}
ds_{h}^{2}=\left( -d\theta ^{2}+\cosh ^{2}\theta f_{AB}y^{A}y^{B}\right) .
\label{LES.02A}
\end{equation}

The reduced equation of (\ref{LES.02}) under the Lie point symmetry
generated by the gradient HV is Laplace equation in the space (\ref{LES.02A}%
). For this reduction we do not have inherited symmetries and there exist
Type II hidden symmetries as stated in section \ref{grHV}.

Before we reduce (\ref{LES.02}) with the symmetry generated from the sp.CKV~$%
X^{3}$, it is best to write the metric (\ref{LES.01}) in new coordinates. We
introduce the new variable $x$ via the relation 
\begin{equation}
z=\sqrt{\frac{R\left( xR-1\right) }{x}}.  \label{LES.03a}
\end{equation}%
In the new variables the Laplacian (\ref{LES.02}) becomes%
\begin{eqnarray}
0 &=&\frac{x^{2}}{R^{2}}u_{xx}-2\frac{x}{R}\left( 2xR-1\right) u_{xR}+u_{RR}+%
\frac{1}{R^{2}}f^{AB}u_{AB}+  \label{LES.04} \\
&&+\frac{\left( m-1\right) }{R}u_{R}-\frac{x}{R^{2}}\left( m-1\right) \left(
2xR-1\right) u_{x}-\frac{1}{R^{2}}\Gamma ^{A}u_{A}  \notag
\end{eqnarray}%
whereas the Lie symmetry $X^{3}$ becomes 
\begin{equation}
X^{3}=\sqrt{\frac{R\left( xR-1\right) }{x}}R\partial _{R}+2p\sqrt{\frac{%
R\left( xR-1\right) }{x}}u\partial _{u}.  \label{LES.04b}
\end{equation}%
The zero order invariants of (\ref{LES.04b}) are $x,y^{A}~,~w=uR^{-2p}.$ We
choose $x,y^{A}$ to be the independent variables and $w=w\left(
x,y^{A}\right) $ the dependent one. Replacing in (\ref{LES.04}) we find the
reduced equation:%
\begin{equation}
x^{2}w_{xx}+f^{AB}w_{AB}-\Gamma ^{A}w_{A}-2p\left( 2p+1\right) w=0
\label{LES.07}
\end{equation}%
We consider cases.\newline
\textbf{The case $m\geq 4.$}

If $2p+1\neq 0,~m\geq 4~$\ then (\ref{LES.07}) becomes%
\begin{equation}
_{\left( m\geq 4\right) }\bar{\Delta}w-2p\left( 2p+1\right) V\left( x\right)
w=0  \label{LES.11}
\end{equation}%
where~$V\left( x\right) =x^{\frac{2}{2-m}}$ and $_{\left( m\geq 4\right) }%
\bar{\Delta}$ is the Laplace operator for the metric%
\begin{equation}
d\bar{s}_{\left( m\geq 4\right) }^{2}=\frac{1}{V\left( x\right) }\left( 
\frac{1}{x^{2}}dx^{2}+f_{AB}dy^{A}dy^{B}\right) .  \label{LES.11b}
\end{equation}

Considering the new transformation $\phi =\int \sqrt{\frac{1}{xV}}dx~\ $~or $%
x=\left( m-2\right) ^{2-m}\phi ^{m-2}$ the metric (\ref{LES.11b}) is written%
\begin{equation}
d\bar{s}_{\left( m\geq 4\right) }^{2}=d\phi ^{2}+\phi ^{2}\bar{f}%
_{AB}dy^{A}dy^{B}  \label{LES.11c}
\end{equation}%
where $\bar{f}_{AB}=\left( m-2\right) ^{-2}f_{AB}$ whereas the potential $%
V\left( \phi \right) =\frac{\left( 2-m\right) ^{2}}{\phi ^{2}}$ which is the
well known Ermakov potential \cite{LeachAndriopoulos}.

This means that the gradient HV $\phi \partial _{\phi }~,~\psi _{\phi }=1$
of the metric (\ref{LES.11c})~satisfies condition (\ref{PE.14}) of Theorem %
\ref{KG} hence it is a Lie point symmetry of (\ref{LES.11}), which is the
Lie symmetry $X^{2}~$ of (\ref{LES.02c}). \newline
Therefore if the metric (\ref{LES.11c}) admits proper CKVs which satisfy the
conditions of Theorem \ref{KG}, then these vectors generate Type II hidden
symmetries of (\ref{LES.02}).\newline
\textbf{The case $m=3.$}

If $2p+1=0~$then $~m=3$ and $f_{AB}$ is a two dimensional metric. In this
case (\ref{LES.07}) becomes%
\begin{equation}
x^{2}w_{xx}+f^{AB}w_{AB}-\Gamma ^{A}w_{A}=0  \label{LES.08}
\end{equation}%
or, by multiplying with ~$x^{2}$%
\begin{equation}
_{\left( m=3\right) }\bar{\Delta}w=0  \label{LES.09}
\end{equation}%
which is the Laplacian in the three dimensional space with metric 
\begin{equation}
d\bar{s}_{\left( m=3\right) }^{2}=\frac{1}{x^{4}}dx^{2}+\frac{1}{x^{2}}%
f_{AB}dy^{A}dy^{B}.  \label{LES.10}
\end{equation}

By making the new transformation $x=\frac{1}{\phi },~\ $ (\ref{LES.09}) is
of the form (\ref{LEH.01}) and admits the gradient HV $\phi \partial _{\phi
}~$which gives an inherited symmetry.

We conclude that Type II hidden symmetries of (\ref{LES.10}) will be
generated from the proper CKVs of the metric (\ref{LES.10}) which satisfy
the conditions of Theorem \ref{BOZKOV}.\newline
\textbf{The case $m=2.$}

\label{spCKVm2}

For $m=2,~$ $f_{AB}$ \ is a one dimensional metric and (\ref{LES.01}) is 
\begin{equation}
ds^{2}=-dz^{2}+dR^{2}+R^{2}d\theta ^{2}
\end{equation}%
which is a flat metric\footnote{%
The only three dimensional space which admits sp.CKV is the flat space,
because in that case we also have a gradient HV\ and a gradient KV.}. In
this space Laplace equation (\ref{LES.02}) admits ten Lie point symmetries,
as many as the elements of the conformal algebra of the flat 3d-space. Six
of these vectors are KVs, one vector is a gradient HV and three are sp.CKVs
(see Appendix). We reduce Laplace equation with the symmetry $X^{3}$~\ and
the reduced equation is (\ref{LES.07}) which for $f_{AB}=\delta _{\theta
\theta }$ becomes 
\begin{equation}
x^{2}w_{xx}+w_{\theta \theta }+\frac{1}{4}w=0.  \label{LES.072}
\end{equation}%
Equation (\ref{LES.072}) is in the form of (\ref{PE.3}) with $%
A^{ij}=diag\left( x^{2},1\right) $ and $B^{i}=0.~$ Replacing in the symmetry
conditions (\ref{PE.3a})-(\ref{PE.3c}) we find the Lie symmetries%
\begin{equation*}
X=\xi ^{i}\partial _{i}+\left( a_{0}w+b\right) \partial _{w}
\end{equation*}%
where $\xi ^{i}$ are the CKVs of the two dimensional space with metric $%
A^{ij}$. In this case all proper CKVs of the two dimensional space $A^{ij}$
generate Type II Lie symmetries. We recall that the conformal algebra of a
two dimensional space is infinite dimensional.

\section{The Wave equation}

\label{Applications}

In this section we apply the general results of the previous section to the
3+1 wave equation in Minkowski spacetime $M^{4}.$ As will be shown we
recover the results of previous studies \cite{AbrahamGovinderArrige2006}
easily and in a straightforward manner. We also amend some of them.

Laplace equation in $M^{4}$ 
\begin{equation}
ds^{2}=-dt^{2}+dx^{2}+dy^{2}+dz^{2}  \label{Ap1.01}
\end{equation}%
is the wave equation in $E^{3}$%
\begin{equation}
u_{tt}-u_{xx}-u_{yy}-u_{zz}=0.  \label{Ap1.02}
\end{equation}

The conformal algebra of the metric (\ref{Ap1.01}) is generated by 15 vector
fields (see Appendix). From Theorem \ref{BOZKOV} we have that the extra Lie
point symmetries of (\ref{Ap1.02}) are the following vector fields 
\begin{equation}
K_{G}^{1}~,~K_{G}^{A}~,~X_{R}^{1A}~,~X_{R}^{AB}~,~H~  \label{Ap1.02A}
\end{equation}%
\begin{equation}
X_{L}^{1}=X_{C}^{1}-tu\partial _{u}~,~X_{L}^{A}=X_{C}^{A}-y^{A}u\partial _{u}
\end{equation}%
where $y^{A}=\left( x,y,z\right) $ with nonzero commutators%
\begin{eqnarray*}
\left[ K_{G}^{I},X_{R}^{IJ}\right] &=&-K_{G}^{J}~,~\left[ K_{G}^{I},~H\right]
=K_{G}^{I} \\
\left[ K_{G}^{I},X_{L}^{A}\right] &=&H-u\partial _{u}~,~\left[
K_{G}^{I},X_{L}^{J}\right] =X_{R}^{IJ} \\
\left[ H,X_{L}^{I}\right] &=&X_{L}^{I}~,~\left[ X_{R}^{IJ},X_{L}^{I}\right]
=X_{L}^{J}~
\end{eqnarray*}%
and the commutators of the rotations $X_{R}^{1A},X_{R}^{AB}$.

\subsection{Reduction with a gradient KV}

\label{M4G}

We consider reduction of (\ref{Ap1.02}) with the gradient KV $%
K_{G}^{z}=\partial _{z}$. The reduced equation is%
\begin{equation}
w_{tt}-w_{xx}-w_{yy}=0  \label{Ap1.033}
\end{equation}%
which is Laplace equation in the space $M^{3}.$ The extra Lie point
symmetries of (\ref{Ap1.033}) are 
\begin{equation}
K_{G}^{1}~,~K_{G}^{x}~,~K_{G}^{y}~,~~X_{R}^{1a}~\ ,~X_{R}^{ab}~,~H
\label{Ap1.033A}
\end{equation}%
and are inherited symmetries (see also the relevant commutators). The Type
II symmetries are the vector fields%
\begin{equation}
\bar{X}_{C}^{1}-\frac{1}{2}tu\partial _{u}~,~\bar{X}_{C}^{x}-\frac{1}{2}%
xu\partial _{u}~,~\bar{X}_{C}^{y}-\frac{1}{2}yu\partial _{u}
\label{Ap1.033B}
\end{equation}%
that is, the Type II hidden symmetries are generated from the sp.CKVs of $%
M^{3}$.

\subsection{Reduction with the gradient HV}

In this case it is better to switch to hyperspherical coordinates $(r,\theta
,\phi ,\zeta ).$ In these coordinates the metric (\ref{Ap1.01}) is 
\begin{equation}
ds^{2}=dr^{2}-r^{2}\left[ d\theta ^{2}+\cosh ^{2}\theta \left( d\phi
^{2}+\cosh ^{2}\phi d\zeta ^{2}\right) \right]  \label{Ap1.04}
\end{equation}%
and the wave equation (\ref{Ap1.02}) becomes%
\begin{equation}
u_{rr}-\frac{1}{r^{2}}\left( u_{\theta \theta }+\frac{u_{\phi \phi }}{\cosh
^{2}\theta }+\frac{u_{\zeta \zeta }}{\cosh ^{2}\theta \cosh \phi ^{2}}%
\right) +\frac{3}{r}u_{r}-2\frac{\tanh \theta }{r^{2}}u_{\theta }-\frac{%
\tanh \phi }{r^{2}\cosh ^{2}\theta }u_{\phi }=0.  \label{Ap1.05}
\end{equation}

According to the analysis of section \ref{grHV} the reduced equation is (\ref%
{LEH.03}), which is the Laplacian in the three dimensional space of the
variables $\left( \theta ,\phi ,\zeta \right) :$%
\begin{equation}
w_{\theta \theta }+\frac{w_{\phi \phi }}{\cosh ^{2}\theta }+\frac{w_{\zeta
\zeta }}{\cosh ^{2}\theta \cosh \phi ^{2}}+2\frac{\tanh \theta }{r^{2}}%
w_{\theta }+\frac{\tanh \phi }{r^{2}\cosh ^{2}\theta }w_{\phi }=0.
\label{Ap1.06}
\end{equation}%
This space is a space of constant curvature. The conformal algebra\footnote{%
All spaces of constant curvature are conformally flat hence they admit the
same conformal algebra with the flat space but in general with different KVs
and CKVs.} of a 3d space of constant non-vanishing curvature consists of 6
non-gradient KVs and 4 proper CKVs \cite{Hallbook,TNA} which in Cartesian
coordinates are given in Appendix. (The rotations and the sp.CKVs of the
flat space are KVs for the space of constant curvature, the rest are proper
gradient CKVs). The conformal factors of the CKVs do not satisfy the
condition $_{h}\Delta \psi =0$ , hence they do not generate Lie point
symmetries for the reduced equation (\ref{Ap1.06}) (see theorem \ref{BOZKOV}%
) whereas for the same reason the KVs are Lie point symmetries of (\ref%
{Ap1.06}). Therefore all point Lie symmetries are inherited and we do not
have Type II\ hidden symmetries.

We note that the proper CKVs of a space of constant non-vanishing curvature
are gradient and their conformal factor satisfies the relation \cite%
{Hallbook} 
\begin{equation*}
\psi _{;ab}=\frac{\psi R}{n-1}g_{ab}\rightarrow g^{ab}\psi _{;ab}=\frac{n}{%
n-1}R\psi \rightarrow _{h}\Delta \psi =\frac{n}{n-1}R\psi
\end{equation*}%
where $R$ is the Ricci scalar of the space of constant curvature. This
implies that they are Lie symmetries of the conformally invariant
Laplace-Beltrami operator \cite{RyanM} but not of the Laplace equation (\ref%
{Ap1.06}).

\subsection{Reduction with a sp.CKV}

Following the steps of section \ref{spCKV} we consider the transformation to
axi-symmetric coordinates $\left( t,R,\theta ,\phi \right) $ in which (\ref%
{Ap1.02}) takes the form 
\begin{equation}
u_{tt}-u_{RR}-\frac{1}{R^{2}}\left( u_{\theta \theta }+\frac{u_{\phi \phi }}{%
\cosh ^{2}\theta }\right) -\frac{2}{R}u_{r}-\frac{\tanh \theta }{R^{2}}%
u_{\theta }=0.  \label{Ap1.07}
\end{equation}%
Applying the transformation (\ref{LES.03a}) $t=\sqrt{\frac{R\left( \tau
R-1\right) }{\tau }}$ we find (note that this is the case $m=3$) that (\ref%
{Ap1.07}) is written as (\ref{LES.04}) and the reduced equation is the
Laplacian $_{\left( m=3\right) }\Delta w~\ \ $for the 3d metric 
\begin{equation}
ds^{2}=\frac{1}{\tau ^{4}}d\tau ^{2}-\frac{1}{\tau ^{2}}\left( d\theta
^{2}+\cosh ^{2}\theta d\phi ^{2}\right) .  \label{Ap.1.08}
\end{equation}%
However the metric (\ref{Ap.1.08}) under the coordinate transformation $\tau
=\frac{1}{T}$ is written%
\begin{equation}
ds^{2}=dT^{2}-T^{2}\left( d\theta ^{2}+\cosh ^{2}\theta d\phi ^{2}\right)
\label{AP1.08}
\end{equation}%
which is the flat 3d Lorentzian metric, which does not admit proper CKVs.
This implies that the Lie point symmetries of the reduced equation are
generated from the KVs/HV/sp.CKVs of the flat $M^{3}$ metric and all are
inherited. Therefore we do not have Type II\ hidden symmetries for the
reduction with a sp.CKV.

The reduction of the 3+1 and 2+1 wave equation has been done previously by
Abraham-Shrauner et. all \cite{AbrahamGovinderArrige2006} and our results
coincide with theirs. We note that in \cite{AbrahamGovinderArrige2006} they
use algebraic computing to find the Lie symmetry generators whereas our
approach is geometric and general and makes no use of algebraic computing
programs. Furthermore our analysis can be generalized to higher dimensions -
where algebraic computing is rather cumbersome - in a straightforward manner.

In the following sections we apply theorem \ref{BOZKOV} to determine the Lie
point symmetries and in addition determine the Type II hidden symmetries of
some spacetimes which are of interest in General Relativity.

\section{LRS Spacetime}

\label{ExLRS}

A spatially Locally Rotational Symmetric (LRS) spacetime has a metric which
admits a group of motions $G_{4}$ acting transitively on spacelike
hypersurfaces $S_{3}$. In coordinates $\{t,x,y,z\}$ the classes of metrics
describing the LRS spacetimes are the following \cite{StephaniEx}:

\begin{equation}
ds^2=\varepsilon [-dt^2+A^2(t)dx^2]+B^2(t)\left[ dy^2+\Sigma ^2(y,k)dz^2%
\right]  \label{sx1.1}
\end{equation}

\begin{equation}
ds^2=\varepsilon \left\{ -dt^2+A^2(t)\left[ dx+\Lambda (k,y)dz\right]
^2\right\} +B^2(t)\left[ dy^2+\Sigma ^2(y,k)dz^2\right]  \label{sx1.2}
\end{equation}

\begin{equation}
ds^{2}=\varepsilon \lbrack
-dt^{2}+A^{2}(t)dx^{2}]+B^{2}(t)e^{2x}(dy^{2}+dz^{2})  \label{sx1.3}
\end{equation}%
where $\varepsilon =\pm 1,\Sigma (y,k)=\sin y,\sinh y,y$ and $\Lambda
(k,y)=\cos y,\cosh y,y$ for $k=1,-1,0$ respectively. (The factor $%
\varepsilon =\pm 1$ essentially distinguishes between the static and the
non-static cases as it can be seen by interchanging the coordinates $t,x$).
These metrics are quite general and contain many well known and important
families of spacetimes. They contain the diagonal Bianchi Type I, III
metrics (with two of the three spacelike metric components equal), the
static spherically/hyperbolically/plane symmetric metrics (interchange $t,x$
in (\ref{sx1.1}) and the signs of $dt$, $dx$), some of the \{2+2\} and
\{1+3\} decomposable metrics etc. The LRS metrics have been classified by
Ellis \cite{I2}. The metrics (\ref{sx1.2}) with $\varepsilon =1$ in Ellis
classification are class $I$ LRS metrics, the metrics (\ref{sx1.1}) and (\ref%
{sx1.3}) are class $II$ LRS and the metrics (\ref{sx1.2}) with $\varepsilon
=-1$ are the class $III$ LRS metrics.

In this section we consider the particular LRS space-time with line element 
\begin{equation}
ds^{2}=-dt^{2}+dR^{2}+R^{s}\left( dz^{2}+dy^{2}\right)  \label{lrs.01}
\end{equation}

When $s=0~$ the space-time (\ref{lrs.01}) reduces to Minkowski space $M^{4}$
considered in the last section.

For $s\neq 0,2 ~$ spacetime (\ref{lrs.01}) admits four KVs \cite{ApolT} 
\begin{equation*}
K^{1}=\partial _{t}~,~K^{2}~,~K^{3}~,~K^{4}
\end{equation*}%
where $K^{1}$ is a gradient KV, $K^{2-4}$ are the elements of $E^{2}$ Lie
algebra, that is,%
\begin{equation*}
K^{2}=\sin y\partial _{z}+\frac{\cos y}{z}\partial _{z}~,~K^{3}=\cos
y\partial _{z}-\frac{\sin y}{z}\partial _{y}~~,~K^{4}=\partial _{y}
\end{equation*}%
and one non-gradient HV 
\begin{equation*}
H=t\partial _{t}+R\partial _{R}+\frac{\left( 2-n\right) }{2}\left( z\partial
_{z}+y\partial _{y}\right) ~,~\psi _{H}=1.
\end{equation*}%
Moreover, in the special case $s=2$ the metric (\ref{lrs.01}) admits the
special CKV%
\begin{equation*}
C_{sp}=\frac{t^{2}+R^{2}}{2}\partial _{t}+tR\partial _{R}~\ ,~\psi _{C}=t
\end{equation*}%
and in that case, $H~$ becomes a gradient HV. In all cases the metric (\ref%
{lrs.01}) does not admit proper non special CKVs.

Laplace equation (\ref{PE.9}) for the line element (\ref{lrs.01}) becomes%
\begin{equation}
-u_{tt}+u_{RR}+\frac{1}{R^{s}}\left( u_{zz}+u_{yy}\right) +\frac{s}{R}%
u_{R}=0.  \label{lrs.04}
\end{equation}

From theorem \ref{BOZKOV} we have that the extra Lie point symmetries of (%
\ref{lrs.04}) are the vector fields%
\begin{equation*}
K^{1}~,~K^{2-4}~,~H\text{ ~for }s\neq 0,2~
\end{equation*}%
\begin{equation*}
K^{1}~,~K^{2-4}~,~H\text{ ,~}~C_{sp}-tu\partial _{u}\text{~for }s=2\text{ }
\end{equation*}%
with non zero commutators:

For $s\neq 0,2$%
\begin{equation*}
\left[ K^{1},H\right] =K^{1}~,~~\left[ K^{2},K^{4}\right] =-K^{3}~,~\left[
K^{3},K^{4}\right] =K^{2}
\end{equation*}%
\begin{equation*}
\left[ K^{2},H\right] =\left( 1-\frac{s}{2}\right) K^{2}\ ,~~\left[ K^{3},H%
\right] =\left( 1-\frac{s}{2}\right) K^{3}~
\end{equation*}%
and for $s=2$%
\begin{equation*}
\left[ K^{1},H\right] =K^{1},~\left[ K^{1},C_{sp}-tu\partial _{u}\right]
=H-u\partial _{u}
\end{equation*}%
\begin{equation*}
\left[ H,C_{sp}-tu\partial _{u}\right] =C_{sp}-tu\partial _{u}~,~\left[
K^{2},K^{4}\right] =-K^{3}~,~\left[ K^{3},K^{4}\right] =K^{2}.
\end{equation*}

Below we study the reduction of (\ref{lrs.04}) using the zero order
invariants generated by the gradient KV $K^{1},$ the HV $H$ and the sp.CKV $%
C_{sp}$.

\subsection{Reduction with a gradient KV}

The gradient KV $K^{1}$ is a non null vector field and the results of
section \ref{grKV} apply. Therefore the reduced equation of (\ref{lrs.04})
is Laplace equation of the 3d space 
\begin{equation}
ds_{\left( 3\right) }^{2}=dR^{2}+R^{s}\left( dz^{2}+dy^{2}\right)
\label{lrs.05a}
\end{equation}%
that is, 
\begin{equation}
w_{RR}+\frac{1}{R^{s}}\left( w_{zz}+w_{yy}\right) +\frac{s}{R}w_{R}=0.
\label{lrs.05}
\end{equation}

The 3d space (\ref{lrs.05a}) is a conformally flat space hence admits a ten
dimensional conformal algebra. Applying theorem \ref{BOZKOV} we find that
the extra Lie symmetries of (\ref{lrs.05}) are:

For $s\neq 0,2$%
\begin{equation*}
w\partial _{w}~\ ,~b\left( \bar{x},z,y\right) \partial _{w}
\end{equation*}%
\begin{equation*}
H_{\left( 3\right) }=R\partial _{R}+\frac{\left( 2-n\right) }{2}\left(
z\partial _{z}+y\partial _{y}\right) ~,~K^{2},~\ K^{3}~,~K^{4}
\end{equation*}%
\begin{equation*}
C^{1}=Rz\partial _{R}+\left[ \frac{2-s}{4}\left( z^{2}-y^{2}\right) -\frac{%
R^{2-s}}{2-s}\right] \partial _{z}+\frac{2-s}{2}zy\partial _{y}-\frac{1}{2}%
zw\partial _{w}
\end{equation*}%
\begin{equation*}
C^{2}=Ry\partial _{R}+\frac{2-s}{2}zy\partial _{z}+\left[ \frac{2-s}{4}%
\left( y^{2}-z^{2}\right) -\frac{R^{2-s}}{2-s}\right] \partial _{y}-\frac{1}{%
2}yw\partial _{w}
\end{equation*}%
and for $s=2~$%
\begin{equation*}
w\partial _{w}~\ ,~b\left( \bar{x},z,y\right) \partial _{w}
\end{equation*}%
\begin{equation*}
H_{\left( 3\right) }=R\partial _{R}~,~K^{2},~\ K^{3}~,~K^{4}
\end{equation*}%
\begin{equation*}
C^{1}=Rz\partial _{R}-\ln R\partial _{z}-\frac{1}{2}zw\partial
_{w}~,~~C^{2}=Ry\partial _{R}-\ln R\partial _{y}-\frac{1}{2}yw\partial _{w}.
\end{equation*}

We observe that the Lie symmetries $H_{\left( 3\right) },~K^{2-4}$ are
inherited symmetries. Lie point symmetries $C^{1-2}$ are generated from the
proper CKVs of (\ref{lrs.05a}); therefore symmetries $C^{1-2}$ are Type II
hidden symmetries.

\subsection{Reduction with a HV}

In this subsection we reduce Laplace equation (\ref{lrs.04}) using the Lie
point symmetry generated by the HV $H$. Recall that $H$ is a gradient HV
only for $s=2$.

Under the coordinate transformation%
\begin{equation*}
t=e^{Sr}\sinh \theta ~~,~~R=e^{r}\cosh \theta
\end{equation*}%
\begin{equation*}
z=\zeta e^{Sr}~~,~y=ve^{Sr}
\end{equation*}%
where $S=\frac{2-n}{2}~$the line element (\ref{lrs.01}) becomes%
\begin{eqnarray*}
ds^{2}&=&-\left( 1-\frac{\left( 2-s\right) ^{2}}{4}\zeta ^{2}\cosh
^{s}\theta \right) e^{2r}dr^{2} \\
&&+e^{2r}d\theta +\left( 2-s\right) e^{2r}\cosh ^{s}\theta drd\zeta
+e^{2r}\cosh ^{s}\theta \left( d\zeta ^{2}+\zeta ^{2}dv^{2}\right)
\end{eqnarray*}%
and Laplace equation is written as follows%
\begin{eqnarray}
0 &=&-u_{rr}+u_{\theta \theta }+\left( 2-s\right) e^{2r}u_{r\zeta }+\frac{1}{%
\cosh ^{s}\theta }\left[ \left( 1-\frac{\left( 2-s\right) ^{2}}{4}\zeta
^{2}\cosh ^{s}\theta \right) u_{\zeta \zeta }+\frac{1}{\zeta ^{2}}u_{vv}%
\right]  \notag \\
&&~~+s\tanh \theta u_{\theta }-su_{R}+\frac{\left( \frac{1}{4}\zeta
^{2}\cosh ^{s}\left( 3s-2\right) \left( 2-s\right) -1\right) }{\zeta \cosh
^{s}\theta }u_{\zeta }  \label{lrs.06}
\end{eqnarray}%
In these coordinates the HV becomes $H=\partial _{r}$. Hence the zero order
invariants are $\alpha ,\beta ,\gamma ~$and $w$. We choose $\alpha ,\beta
,\gamma $ as the independent variables and $w=w\left( \alpha ,\beta ,\gamma
\right) $ as the dependent variable and find the reduced equation%
\begin{eqnarray}
0 &=&w_{\alpha \alpha }+\frac{1}{\cosh ^{s}\alpha }\left[ \left( 1-\frac{%
\left( 2-s\right) ^{2}}{4}\beta ^{2}\cosh ^{s}\alpha \right) w_{\beta \beta
}+\frac{1}{\beta ^{2}}w_{\gamma \gamma }\right]  \notag \\
~~ &&~~~~~~~~~~~~~~+s\tanh \alpha w_{\alpha }+\frac{\left( \frac{1}{4}\beta
^{2}\cosh ^{s}\left( 3s-2\right) \left( 2-s\right) -1\right) }{\beta \cosh
^{s}\alpha }w_{\beta }  \label{lrs.07}
\end{eqnarray}

Equation (\ref{lrs.07}) admits the Lie point symmetry $K^{4}$ for $s\neq
0,2~ $\ plus the symmetries $w\partial _{w},~b\left( \alpha ,\beta ,\gamma
\right) \partial _{w}$ because (\ref{lrs.07}) is linear in $w$. $K^{4}$ is
inherited therefore we do not have Type II hidden symmetries.

In the case where $s=2$, $H$ is a gradient HV and the results of section \ref%
{grHV} apply. In that case (\ref{lrs.07}) becomes%
\begin{equation}
_{h}\bar{\Delta}w=0  \label{lrs.08}
\end{equation}%
where $_{h}\bar{\Delta}$ is Laplace operator for the 3d metric%
\begin{equation}
d\bar{s}_{\left( 3\right) }^{2}=d\alpha ^{2}+\cosh ^{2}\alpha \left( d\beta
^{2}+\beta ^{2}d\gamma ^{2}\right) .  \label{lrs.09}
\end{equation}%
This metric is conformally flat hence admits a ten dimensional conformal
algebra with a three dimensional Killing algebra. From theorem \ref{BOZKOV},
we have that (\ref{lrs.08}) admits as Lie point symmetries the $E^{2}$ Lie
algebra plus the vectors $w\partial _{w},~b\left( \alpha ,\beta ,\gamma
\right) \partial _{w}$. In that case all symmetries are inherited and we do
not have Type II hidden symmetries.

\subsection{Reduction with a sp.CKV}

In the case $s=2$ Laplace equation admits a Lie point symmetry generated by
the sp.CKV $C_{sp}$ and the results of section \ref{spCKV} apply. That is,
from subsection \ref{spCKV} and for $m=3$ the reduced equation is 
\begin{equation}
_{\left( sp\right) }\Delta w=0  \label{lrs.10}
\end{equation}%
where $_{\left( sp\right) }\Delta $ is the Laplace operator for the metric~ 
\begin{equation}
d\bar{s}_{\left( m=3\right) }^{2}=d\bar{x}^{2}+\bar{x}^{2}\left(
dz^{2}+dy^{2}\right) .  \label{lrs.11}
\end{equation}%
and we have defined $x=\frac{1}{\bar{x}}$. From theorem \ref{BOZKOV} we have
that (\ref{lrs.10}) admits the Lie point symmetries~%
\begin{equation*}
\bar{H}_{\left( 3\right) }=\partial _{\bar{x}}~~,~~w\partial _{w}~\
,~b\left( \bar{x},z,y\right) \partial _{w}
\end{equation*}%
\begin{equation*}
\bar{C}^{1}=\bar{x}z\partial _{\bar{x}}-\ln \bar{x}\partial _{z}-\frac{1}{2}%
zw\partial _{w}~,~~\bar{C}^{2}=\bar{x}y\partial _{\bar{x}}-\ln \bar{x}%
\partial _{y}-\frac{1}{2}yw\partial _{w}.
\end{equation*}

$H_{\left( 3\right) }$ is a gradient HV for space (\ref{lrs.11}) and the Lie
point symmetries $\bar{C}^{1-2}$ are generated by proper CKVs of (\ref%
{lrs.11}). From subsection \ref{spCKV} we conclude have that $\bar{H}%
_{\left( 3\right) }$ is an inherited symmetry whereas $\bar{C}^{1-2}$ are
Type II hidden symmetries.

\section{The algebraically special empty space Petrov Type III solution}

\label{The Heat equation in spaces which admit a non-gradient HV}

In this section we consider the reduction of Laplace equation in spaces
which do not admit gradient KVs or a gradient HV. As it has been mentioned
in section \ref{Laplace equation in certain Riemanian spaces} we shall
consider the algebraically special solutions of Einstein equations which
admit a homothetic algebra acting simply transitively. These spacetimes have
been determined in \cite{Steele1991(b)} and are of Petrov type D,N,II and
III In the following we restrict our discussion to Petrov Type III only with
metric (\ref{P3.00}), because both the method of work and the results are
the same for the remaining Petrov types in this class of spacetimes.
Spacetime (\ref{P3.00}) admits the four dimensional conformal algebra
generated by the vector fields%
\begin{eqnarray*}
K^{1} &=&\partial _{\rho }~,~K^{2}=\partial _{y}~,~K^{3}=v\partial _{v}-\rho
\partial _{\rho }+2x\partial _{x}+2y\partial _{y} \\
H &=&v\partial _{v}+\rho \partial _{\rho }~,~\psi =1
\end{eqnarray*}%
where $K^{1-3}$ are KVs and $H$ is a non-gradient HV. (The space does not
admit proper CKVs).

In this spacetime the Laplacian takes the form:%
\begin{equation}
-\frac{3}{2}xu_{vv}+2u_{v\rho }+\frac{x^{3}}{v^{2}}\left(
u_{xx}+u_{yy}\right) -3\frac{x}{v}u_{v}+\frac{2}{v}u_{\rho }=0.
\label{P3.01}
\end{equation}%
From{\LARGE \ }theorem{\LARGE \ }\ref{BOZKOV} we have{\LARGE \ }that the
extra Lie point symmetries are the vector fields 
\begin{equation*}
X_{1-3}=K_{1-3}~,~~X_{4}=H
\end{equation*}%
with nonzero commutators:%
\begin{equation*}
\left[ X_{2,},X_{3}\right] =2X_{2}
\end{equation*}%
\begin{equation*}
\left[ X_{3},X_{1}\right] =X_{1}~,~~\left[ X_{1},X_{4}\right] =X_{1}.~
\end{equation*}%
We use $X_{4}~$to reduce the PDE because this is the Lie symmetry generated
by the nongradient HV{\LARGE .}

The zero order invariants of $X_{4}$ are $\sigma =\frac{\rho }{v},x,y,w$. We
choose $\sigma ,x,y$ as the independent variables and $w=w\left( \sigma
,x,y\right) $ as the dependent variable and we find the reduced equation 
\begin{equation}
-\sigma \left( \frac{3}{2}x\sigma +2\right) w_{\sigma \sigma }+x^{3}\left(
w_{xx}+w_{yy}\right) =0.  \label{P3.02}
\end{equation}%
Equation (\ref{P3.02}) can be written 
\begin{equation}
_{III}\Delta ^{\ast }w-\left( \frac{3x\sigma }{2}+1\right) w_{\sigma }-\frac{%
3x^{3}\sigma }{2\left( 3x\sigma +4\right) }w_{x}=0  \label{P3.03}
\end{equation}%
where $_{III}\Delta ^{\ast }$ is the Laplacian for the metric%
\begin{equation}
ds^{2}=-\frac{1}{\sigma \left( \frac{3}{2}x\sigma +2\right) }d\sigma ^{2}+%
\frac{1}{x^{3}}\left( dx^{2}+dy^{2}\right) .  \label{P3.04}
\end{equation}%
The Lie point symmetries of (\ref{P3.03}) are generated from the conformal
algebra of (\ref{P3.04}) with some extra conditions (see eqs. (\ref{PE.3a})-(%
\ref{PE.3c})). We find that equation (\ref{P3.03}) admits as Lie point
symmetries the vectors~$\partial _{y}~,~x\partial _{x}+y\partial _{x}-\sigma
\partial _{\sigma }$ which are inherited symmetries. Hence we do not have
Type II hidden symmetries with this reduction.

\section{The spatially flat $n-$ dimensional FRW like metric}

\label{FRWCKV}

As a final example we consider the $n$ dimensional FRW like space $\left(
n>2\right) $ with metric\footnote{%
This is a special form of the FRW metric which admits a HV.} (\ref{EC.01})
where $\delta _{AB}$ is the $n-1$ dimensional Euclidian metric. The
reduction of Laplace equation in this space (for $n=4$) has been studied
previously in \cite{Kara}. The metric (\ref{EC.01}) is conformally flat
hence admits the same CKVs with the flat space but with different conformal
factors. More precisely the space admits

a. $\left( n-1\right) +\frac{\left( n-2\right) \left( n-3\right) }{2}~$KVs
the $K_{G}^{A},~X_{R}^{AB}$

b. $1$ gradient HV$~$the $K_{G}^{1}=\partial _{t}~~$\newline
the rest vectors being proper\ CKVs \cite{Maartens}. In this space, Laplace
equation (\ref{PE.9}) becomes%
\begin{equation}
e^{-2t}\left[ u_{tt}-\left( \delta ^{AB}u_{AB}\right) +\left( n-2\right)
u_{t}\right] =0  \label{EC.02}
\end{equation}%
and the extra Lie point symmetries are 
\begin{equation*}
K_{G}^{A},~X_{R}^{AB}~,~K_{G}^{1}~
\end{equation*}%
\begin{equation*}
C^{A}=X_{R}^{1A}-2py^{A}u\partial _{u}
\end{equation*}%
where $2p=\frac{2-n}{2}$ and with nonzero commutators%
\begin{equation*}
\left[ K_{G}^{1},C^{A}\right] =K_{G}^{A}~~
\end{equation*}%
\begin{equation*}
\left[ C^{A},C^{B}\right] =-X_{R}^{AB}~,~\left[ X_{R}^{AB},C^{A}\right]
=C^{B}
\end{equation*}%
and the commutators of the rotations $X_{R}^{AB}$.

\subsection{Reduction with the gradient HV}

The gradient HV $K_{G}^{1}=\partial _{t}$ is a Lie point symmetry of the
Laplacian (\ref{EC.02}) hence we consider the reduction by this vector. The
zero order invariants are $y^{A},w$ and lead to the reduced equation 
\begin{equation}
\delta ^{AB}u_{AB}=0  \label{EC.07}
\end{equation}%
which is Laplace equation in the flat space $E^{n-1}$. We consider again
cases.\newline
{\textbf{Case $n>3$}}

In this case the Lie point symmetries of (\ref{EC.07}) are given by the
vectors (see Appendix \ref{CKVsFlat}) 
\begin{equation}
K_{G}^{A}~,~X_{R}^{AB}~,~_{n-1}H~,X_{C}^{A}-y^{A}u\partial _{u}.
\end{equation}%
From these the $K_{G}^{A},X_{R}^{AB}$ are inherited symmetries and the rest
- which are produced by the HV and the sp.CKVs of the space $E^{n-1}$ - are
Type II hidden symmetries.\newline
\textbf{Case} $\mathbf{n=3}$

In this case, the reduced equation (\ref{EC.07}) is the Laplacian in $E^{2},$
hence admits an infinite dimensional Lie algebra \cite{Stephani}. Type II
hidden symmetries are generated from the HV and the CKVs of $E^{2}$.

\subsection{Reduction with a proper CKV}

We consider next the reduction with a proper CKV. We may take any of the
vectors $X_{R}^{1A}$ (because, as one can see in the Appendix, there is a
symmetry between the coordinates $y^{A}$). We choose the vector field 
\begin{equation*}
X_{R}^{1x}=x\partial _{t}+t\partial _{x}+2pxu\partial _{u}.
\end{equation*}%
whose zero order invariants are 
\begin{equation*}
R=t^{2}-x^{2}~,~y^{C}~,~w=e^{-2pt}u.
\end{equation*}%
We take the dependent variable to be the $w=w\left( R,y^{C}\right) $ and
find the reduced equation 
\begin{equation}
4Rw_{RR}-\delta ^{ab}w_{ab}+4w_{R}-4p^{2}w=0  \label{EC.03}
\end{equation}%
where $a=1,\ldots ,n-2$. We consider cases.

{\textbf{Case $n>3.$}}

For $n>3$ equation (\ref{EC.03}) is 
\begin{equation}
_{C}\Delta w-4p^{2}f\left( R\right) w=0  \label{EC.04}
\end{equation}%
where $_{C}\Delta $ is the Laplace operator for the $\left( n-1\right) $
dimensional metric%
\begin{equation}
ds_{C}^{2}=\frac{1}{f\left( R\right) }\left( \frac{1}{4R}dR^{2}-\delta
_{ab}dy^{a}dy^{b}\right)  \label{EC.05}
\end{equation}%
and $f\left( R\right) =R^{-\frac{1}{n-3}}$. The metric (\ref{EC.05}) is
conformally flat hence we know its conformal algebra. Application of theorem %
\ref{KG} gives that the Lie point symmetries of (\ref{EC.04}) are the vector
fields 
\begin{eqnarray*}
X_{u} &=&u\partial _{u}~,~X_{b}=b\partial _{u} \\
X_{K}^{a} &=&\partial _{y^{a}}~,~X_{R}^{ab}~=y^{b}\partial
_{a}-y^{a}\partial _{b}.
\end{eqnarray*}%
These are inherited symmetries (this result agrees with the commutators). We
conclude that for this reduction we do not have Type II hidden symmetries.

{\textbf{Case $n=3.$}}

For $n=3$ the reduced equation is a two dimensional equation $\left( \text{%
that is }\delta _{AB}=\delta _{yy}\right) $%
\begin{equation}
4Rw_{RR}-w_{yy}+4w_{R}-\frac{1}{4}w=0
\end{equation}%
and admits as Lie point symmetry the KV $\partial _{y}$ which is an
inherited symmetry. Hence we do not have Type II hidden symmetries.

We conclude that the reduction of Laplace equation in a $n$ dimensional FRW
like space with the proper CKV does not produce Type II hidden symmetries
and in fact the inherited symmetries of the reduced equation are the KVs of
the flat metric.

\section{Conclusion}

\label{Conclusion}

Up to now the study of Type II\ hidden symmetries has been done by counter
examples or by considering very special PDEs and in low dimensional flat
spaces. In this paper we improve this scenario and study the problem of Type
II\ hidden symmetries of second order PDEs from a geometric point of view in 
$n-$ dimensional Riemannian spaces. We have considered the reduction of
Laplace equation and the consequent possibility of existence of Type II
hidden symmetries in some general classes of spaces which admit some kind of
symmetry hence admit nontrivial Lie symmetries. The general conclusion of
this study is that the Type II\ hidden symmetries of Laplace equation are
directly related to the transition of the CKVs from the space where the
original equation is defined to the space where the reduced equation resides.

We may summarize the general conclusions of this study as follows:

\begin{itemize}
\item If we reduce Laplace equation with a non null gradient KV the reduced
equation is again Laplace equation in the non-decomposable space. In this
case the Type II hidden symmetries are generated from the special and the
proper CKVs of the non-decomposable space.

\item If we reduce Laplace equation with a gradient HV the reduced equation
is again Laplace equation for an appropriate metric. In this case the Type
II hidden symmetries are generated from the HV and the special/proper CKVs.

\item If we reduce Laplace equation with the symmetry generated by a sp.CKV
in a space which admits a non null KV, the reduced equation is the Klein
Gordon equation (\ref{PE.10}) for an appropriate metric which inherits the
Lie symmetry generated by the gradient HV. In this case the Type II
symmetries are generated from the proper CKVs.
\end{itemize}

\textbf{Acknowledgments.} We would like to thank the referee for helpful
comments which have improved the manuscript.

\appendix

\section{Appendix: Conformal algebra of a flat space of Lorentzian /
Euclidian character}

\label{CKVsFlat}

We consider a flat space of dimension $n>2~$ with metric 
\begin{equation*}
ds^{2}=\varepsilon dt^{2}+\delta _{AB}dy^{A}dy^{B}~,~\varepsilon =\pm 1.
\end{equation*}%
The conformal algebra of the space consists of the following vectors\newline

$n~$ gradient KVs 
\begin{equation*}
K_{G}^{1}=\partial _{t}~,~K_{G}^{A}=\partial _{A}
\end{equation*}

$\frac{n\left( n-1\right) }{2}~$ non~gradient KVs (rotations) 
\begin{equation*}
X_{R}^{1A}=y^{A}\partial _{t}-\varepsilon t\partial
_{A}~,~X_{R}^{AB}=y^{B}\partial _{A}-y^{A}\partial _{B}
\end{equation*}

$1~$ gradient HV$~$%
\begin{equation*}
H=t\partial _{t}+\sum\limits_{A}y^{A}\partial _{A}
\end{equation*}

$n~~$sp.CKVs%
\begin{equation*}
X_{C}^{1}=\frac{1}{2}\left( t^{2}-\varepsilon \sum\limits_{A}\left(
y^{A}\right) ^{2}\right) \partial _{t}+t\sum\limits_{A}y^{A}\partial _{A}
\end{equation*}%
\begin{equation*}
X_{C}^{A}=ty^{A}\partial _{t}+\frac{1}{2}\left( \varepsilon t^{2}+\left(
y^{A}\right) ^{2}-\sum\limits_{B\neq A}\left( y^{B}\right) ^{2}\right)
\partial _{A}+y^{A}\sum\limits_{B\neq A}y^{B}\partial _{B}
\end{equation*}%
where $y^{A}=1...n-1~$with conformal factor $\psi _{C}^{1}=t~$and $\psi
_{C}^{A}=y^{A}.~$For $n>2$ the flat space does not admit proper CKVs \cite%
{TNA}.

For $n=2$ the vector field%
\begin{equation*}
X=\left[ f\left( t+i\sqrt{\varepsilon }x\right) -g\left( t-i\sqrt{%
\varepsilon }x\right) +c_{0}\right] \partial _{t}+i\sqrt{\varepsilon }\left[
f\left( t+i\sqrt{\varepsilon }x\right) +g\left( t-i\sqrt{\varepsilon }%
x\right) \right] \partial _{x}
\end{equation*}%
is the generic CKV, that is, includes the KVs, the HV, and the sp.CKVs.

\end{document}